\documentclass[psamsfonts,12pt]{article}
\usepackage[mathscr]{eucal}
\usepackage{graphicx}
\usepackage{amsthm}
\usepackage{amsmath}
\usepackage{latexsym}
\usepackage{amscd}
\usepackage{amsfonts}
\usepackage{amssymb}
\usepackage{graphicx}
\usepackage{graphicx}
\date{}
\RequirePackage{fix-cm}
\usepackage{graphicx}
\usepackage{amssymb}
\usepackage[mathscr]{eucal}
\usepackage{amsmath}
\usepackage{mathtools}

 \usepackage{mathptmx}      % use Times fonts if available on your TeX system
\usepackage{latexsym}
\usepackage{booktabs}
\usepackage{multirow}
\usepackage{tabularx}
\usepackage{caption}
\usepackage[figuresright]{rotating}
\usepackage{latexsym}

 \usepackage{mathptmx}      % use Times fonts if available on your TeX system
 \usepackage{latexsym}
 \usepackage{booktabs}
 \usepackage{multirow}
 \usepackage{tabularx}
 \usepackage{caption}
 \usepackage[figuresright]{rotating}
 
 \usepackage{graphicx}
 \usepackage[mathscr]{eucal}
 \usepackage{amsmath,amssymb,latexsym,amscd}
 \usepackage{multirow}
 
 \usepackage{mathptmx} 
 \usepackage{latexsym}
 \usepackage{booktabs}
 \usepackage{multirow}
 \usepackage{tabularx}
 \usepackage{caption}
 \usepackage[figuresright]{rotating}

\RequirePackage{fix-cm}
\usepackage{graphicx}
\usepackage{float}

\usepackage{amssymb}
\usepackage[mathscr]{eucal}

\def\C{\mathbb C}

\mathchardef\ssneq="3B28
\def\m@th{\mathsurround=\z@}

\begin{document}

%\title{El Operador de Black-Scholes, $C_{0}$-Semigrupos y Diferentes tipos de Opciones de Inversi\'on}
\title{Solution for the nonlinear relativistic harmonic oscillator via Laplace-Adomian decomposition method}
\author{O. Gonz\'alez-Gaxiola$^{a}$\footnote{ogonzalez@correo.cua.uam.mx},  J. A. Santiago$^{a}$, J. Ruiz de Ch\'avez$^{b}$\\$^{a}$￼Departamento de Matem\'aticas Aplicadas y Sistemas\\ Universidad Aut\'onoma Metropolitana-Cuajimalpa,\\C.P. 05300 Mexico, D.F., Mexico\\ $^{b}$Departamento de Matem\'aticas,\\ Universidad Aut\'onoma Metropolitana-Iztapalapa. \\ San Rafael Atlixco 186, A.P. 55534, Col. Vicentina,\\ Iztapalapa, 09340, M\'exico D.F.} \maketitle

%\addtocontents{toc}{chapter}{\textbf{Introducci\'on}}
%\addcontentsline{toc}{chapter}{Introducci\'on} \textbf{\large
%Introducci\'on} \markboth{Introducci\'on}{Introducci\'on}

\begin{abstract}
 
\noindent Far as we know there are not exact solutions to the equation of motion for
a relativistic harmonic oscillator.	In this paper, the relativistic harmonic oscillator equation which is a nonlinear ordinary differential equation is studied by means of a
combined use of the  Adomian Decomposition Method and the Laplace Transform (LADM). The results that we have obtained, a series of powers of functions, have never been reported
and show a very good match when  compared with other approximate solutions, obtained by different methods. The  method here proposed works with high degree of accuracy and 
because  it requires less computational effort, it is very convenient to solve this kind of nonlinear differential equations.

\end{abstract}

%This poster is partially based on \cite{Osw}.

    %%% Introduction
 %   \vspace{0.8cm}
  %  \begin{center}
       % \pbox{0.8\textwidth}{}%%
      %  {linewidth=2mm,framearc=0.1,linecolor=lightblue,fillstyle=gradient,gradangle=0,%%
     %   gradbegin=white,gradend=whiteblue,gradmidpoint=1.0,framesep=1em}{%%
    %    \begin{center}
         %  \textbf{ Introduction (Mathematical Model)}
   %     \end{center}}
  %  \end{center}
 %   \vspace{1cm}

\section{Introduction}
\label{Intro}
\noindent Many of the phenomena that arise in the real world can be  described by means of nonlinear  partial and ordinary differential equations  and, in some cases,  by integral or differo-integral equations. However, most of the mathematical methods developed so far,  are only capable to  solve linear differential equations.  In the 1980's, George Adomian (1923-1996) introduced a powerful method to solve nonlinear differential equations. Since then, this method is known as the Adomian decomposition method (ADM) \cite{ADM-0,ADM-1}. The technique is based on a decomposition of a solution of a nonlinear differential equation in a series of functions. Each term of the series is obtained from a polynomial generated by a power series expansion of an analytic function. The Adomian method is very simple in an abstract formulation but the difficulty arises in calculating the polynomials that becomes a non-trivial task. This method has widely been used to solve equations that come from nonlinear models as well as to solve fractional differential equations \cite{Das-1,Das-2,Das-3}. The advantage of this method is that, it solves the problem directly without the need of linearization, perturbation, or any other transformation, and also, reduces the massive computation works required by most other methods.

\noindent The relativistic nonlinear harmonic oscillator, was  studied by first time in the middle of the last century \cite{Pen,Gold}. In spite of its importance in several models of physics, 
exact solutions of its equation of motion have not been obtained. 
In the present work we will use the Adomian  decomposition method in combination with the Laplace transform (LADM) \cite{Waz-Lap} to determine the  relativistic oscillator 
solutions. This equation is a nonlinear ordinary differential equation that, in physics, is used to model a simple one dimensional harmonic oscillator  with relativistic velocities.
 We will decomposed the nonlinear terms of this equation using the Adomian polynomials and then, in combination with the use of the Laplace transform, we will obtain an algorithm to solve the problem subject to initial conditions. Finally, we will illustrate our procedure and the quality of the obtained algorithm by solving two examples in which the nonlinear differential equation is solved for  different initial conditions .

\noindent Our work is divided in several sections. In ``The Adomian Decomposition Method Combined With Laplace Transform'' section, we present, in a brief and self-contained manner, the LADM. Several references are given to delve deeper into the subject and to study its mathematical foundation that is beyond the scope of the present work. In ``The Relativistic Harmonic Oscillator'' section, we  also give a brief introduction to the model described by the relativistic harmonic oscillator and we compare our results with the previous one obtained in this respect.
In ``Solution of the Relativistic Harmonic Oscillator Equation Through LADM'' section,
we will establish that LADM can be used to solve this equation in a very simple way. In ``Application to the Relativistic Harmonic Oscillator'' section, we will show by means of two examples, the quality and precision of our method, comparing the obtained results with existing approximate solutions available in the literature and obtained by other methods.
Finally, in the ``Conclusion and Summary"  section, we present our conclusions.

\section{The Adomian Decomposition Method Combined with Laplace Transform}
\label{H-O-0}
\noindent  The ADM is a method  to solve ordinary and nonlinear  differential equations. Using this method is possible to express analytic solutions in terms of a series \cite{ADM-1}. 
In a nutshell,  the method identifies and separates the linear and nonlinear parts of a differential equation. Inverting and applying the highest order differential operator that is contained in the linear part of the equation, it  is possible to express the solution in terms of the rest of the equation affected by the inverse operator.  At this point, the solution is proposed by means of  a series 
with terms that will be  determined and that give rise to the Adomian Polynomials \cite{Waz-0}. The nonlinear part can also be expressed in terms of these polynomials. 
The initial (or the border conditions) and the terms that contain the independent variables will be considered as the initial  approximation. In this way and by means of a  recurrence relations, it is possible to find the  terms of the series that give the approximate solution of the differential equation. In the next paragraph we will see how to use the Adomian  decomposition method in combination with the Laplace transform (LADM).\\

\noindent Let us consider the homogeneous differential equation of second order:
\begin{equation}
\frac{d^{2}x}{dt^2}+N(x)=0\label{eq:y1}
\end{equation}
with initial conditions 
\begin{equation}
x(0)=\alpha,\quad x'(0)=\beta\label{eq:y2}
\end{equation}
where $\alpha$, $\beta$ are real constants and $N$ is a nonlinear operator acting on the dependent variable $x$ and some of its derivatives.\\  
In general, if we consider the second-order differential operator $L_{tt}=\frac{\partial^2}{\partial t^2}$, then the equation (\ref{eq:y1}) could be  written as 
\begin{equation}
L_{tt}x(t)+N(x(t))=0.\label{eq:y3}
\end{equation}
Solving for $L_{tt}x(t)$, we have 
\begin{equation}
L_{tt}x(t)=-N(x(t))\label{eq:y4}.
\end{equation}
The LADM consists of applying Laplace transform (denoted throughout this paper by $\mathcal{L}$) first on both sides of Eq. (\ref{eq:y4}), obtaining
\begin{equation}
\mathcal{L}\{L_{tt}x(t)\}= -\mathcal{L}\{N(x(t))\}.\label{eq:y5}
\end{equation}
An equivalent expression to  (\ref{eq:y5}) is
\begin{equation}
s^{2}x(s)-sx(0)-x'(0)= -\mathcal{L}\{Nx(t)\},\label{eq:y6}
\end{equation}
using the initial conditions (\ref{eq:y2}), we have
\begin{equation}
x(s)=\frac{\alpha}{s}+\frac{\beta}{s^2}-\frac{1}{s^2}\mathcal{L}\{N(x(t))\}\label{eq:y7}
\end{equation}
now, applying the inverse Laplace transform to equation (\ref{eq:y7})
\begin{equation}
x(t)=\alpha+\beta t-\mathcal{L}^{-1}\big[\frac{1}{s^2}\mathcal{L}\{N(x(t))\}\big]. \label{eq:y8}
\end{equation}
\noindent The ADM method proposes a series solution   $x(t)$ given by,
\begin{equation}
x(t)= \sum_{n=0}^{\infty}x_{n}(t).\label{eq:y7-1}
\end{equation}
The nonlinear term $N(x)$ is given by
\begin{equation}
N(x)= \sum_{n=0}^{\infty}A_{n}(x_{0},x_{1},\ldots, x_{n})\label{eq:y8-1}
\end{equation}
where   $\{A_{n}\}_{n=0}^{\infty}$ is the so-called Adomian polynomials sequence established in \cite{Waz-0} and \cite{Ba} and, in general, give us term to term:\\
$A_{0}=N(x_0)$\\
$A_{1}=x_{1}N'(x_0)$\\
$A_{2}=x_{2}N'(x_{0})+\frac{1}{2}x_{1}^{2}N''(x_0)$\\
$A_{3}=x_{3}N'(x_{0})+x_{1}x_{2}N''(x_0)+\frac{1}{3!}x_{1}^{3}N^{(3)}(x_{0})$\\
$A_{4}=x_{4}N'(x_{0})+(\frac{1}{2}x_{2}^{2}+x_{1}x_{3})N''(x_0)+\frac{1}{2!}x_{1}^{2}x_{2}N^{(3)}(x_{0})+\frac{1}{4!}x_{1}^{4}N^{(4)}(x_0)$\\
$ \vdots$.\\
\noindent Other polynomials can be generated in a similar way. Some other approaches to obtain Adomian's polynomials can be found in  \cite{Duan,Duan1}.\\
\noindent Using  (\ref{eq:y7-1}) and (\ref{eq:y8-1}) into equation (\ref{eq:y8}), we obtain,
\begin{equation}
\sum_{n=0}^{\infty}x_{n}(t)= \alpha+\beta t-\mathcal{L}^{-1}\Big[\frac{1}{s^2}\mathcal{L}\{\sum_{n=0}^{\infty}A_{n}(x_{0},x_{1},\ldots, x_{n})\}\Big].\label{eq:y10}
\end{equation}
From the equation (\ref{eq:y10}) we deduce the  recurrence formula:
\begin{equation}
\left\{
    \begin{array}{ll}
x_{0}(t)=\alpha+\beta t,\\
x_{n+1}(t)=-\mathcal{L}^{-1}\Big[\frac{1}{s^2}\mathcal{L}\{A_{n}(x_{0},x_{1},\ldots, x_{n})\}\Big],\;\; n=0,1,2,\ldots
    \end{array}
\right.\label{eq:y11}
\end{equation}
Using  (\ref{eq:y11}) we can obtain an approximate solution of (\ref{eq:y1}), (\ref{eq:y2}) using  
\begin{equation}
x(t)\approx \sum_{n=0}^{k}x_{n}(t),\;\; \mbox{where} \;\; \lim_{k\to\infty}\sum_{n=0}^{k}x_{n}(t)=x(t).\label{eq:y12}
\end{equation}
 It becomes clear that, the Adomian decomposition method, combined with the Laplace transform needs less work in comparison with  the traditional  Adomian decomposition method.  This method decreases considerably the volume of calculations. The decomposition procedure of Adomian will be  easily set, without linearising the problem. With this approach, the solution is found in the form of a convergent series with easily computed components; in many cases, the convergence of this series is very fast and only a few terms are needed in order to have an idea of how the solutions behave. Convergence conditions of this series are examined by several authors,  mainly in \cite{Y3,Y4,Y1,Y2}. Additional references related to the use of the Adomian Decomposition Method, combined with the Laplace transform, can be found in \cite{Waz-Lap,Khu,Y} and references therein.

% by implementing the decomposition method rather than the standard methods for the exact solutions. 

\section{The Relativistic Harmonic Oscillator}
\label{H-O-1}

\noindent The equation of motion of the relativistic harmonic oscillator is given by  the nonlinear differential equation \cite{Bia}:
\begin{equation}\label{Osc-1}
\frac{d^2x}{dt^2}+\Big[1-\Big(\frac{dx}{dt}\Big)^2\Big]^{\frac{3}{2}}x=0,\quad x(0)=0,\quad \frac{dx}{dt}(0)=\beta.
\end{equation}
This normalized, dimensionless form of the equation is based on taking the rest mass $m$ to be unity and the speed of light $c$ to also be unity \cite{Mic1}. It is easy to verify that the dimensionless length $x$ and the dimensionless time $t$ are related to the dimensional variables $\bar{x}$ and $\bar{t}$ through $x=\omega_{0}\bar{x}/c$ and $t=\omega_{0}\bar{t}$, respectively, where $\omega_{0}=\sqrt{k/m}$ is the angular frequency for the non-relativistic oscillator.\\
\noindent As far as we know, no exact solution of the nonlinear equation (\ref{Osc-1}) has yet been published and therefore the research work about equation
 (\ref{Osc-1})  has been intense; a fundamental result reported in \cite{Mic1} 
 is that all the solutions of  (\ref{Osc-1}) are periodic functions with the period dependent of the initial velocity
 $\beta$.
 In the same work, an approximation solution of (\ref{Osc-1})  was found using the harmonic balance method (HBM), it is given by 
 \begin{equation} \label{hbm}
\begin{split}
x_{\mbox{\tiny HBM}}(t) & = \frac{\beta}{\omega}\Big(\frac{3\beta^{4}+8\beta^{2}+64}{64}\Big)\sin(\omega t)-\frac{\beta^3}{24\omega}\Big(\frac{3\beta^{2}+128}{128}\Big)\sin(3\omega t) \\
& +\Big(\frac{3\beta^{5}}{640\omega}\Big)\sin(5\omega t),\;\; \mbox{where}\;\; \omega=\sqrt[4]{\frac{2-2\beta^2}{2-\beta^2}}\;\; \mbox{and}\;\; 0<\beta<1.
\end{split}
\end{equation}
Some more detailed work in the same direction was reported ten years later in \cite{Bel-1,Bel-2}. 
After that, in \cite{Eb}, using the  differential transformation method (DTM), some periodic solutions were obtained and more recently the relativistic harmonic oscillator is
studied by using the homotopy perturbation method (HPM) \cite{Bia}, where a good approximation is obtained using the fact that the solutions are periodic functions.\\
In the following section we will develop an algorithm using the method described in  ``The
Adomian Decomposition Method Combined with Laplace Transform'' section in order to
solve the nonlinear differential equation (\ref{Osc-1}) without resort to any truncation or linearization and not assuming  {\it a priori} that the solutions are periodic functions. 

\section{\bf Solution of the Relativistic Harmonic Oscillator Equation Through LADM}
\label{NLC-ADM}
\noindent Comparing (\ref{Osc-1}) with equation  (\ref{eq:y4}) we have that $L_{tt}$ and $N$ becomes:
\begin{equation}
L_{tt}x=\frac{d^2}{dt^2}x,\;\; Nx=\Big[1-\Big(\frac{dx}{dt}\Big)^2\Big]^{\frac{3}{2}}x .\;\; \label{Oper-1}
\end{equation}
\noindent By using now equation (\ref{eq:y11}) through the  LADM  method we obtain recursively
\begin{equation}
\left\{
    \begin{array}{ll}
x_{0}(t)=\beta t,\\
x_{n+1}(t)=-\mathcal{L}^{-1}\Big[\frac{1}{s^2}\mathcal{L}\{A_{n}(x_{0},x_{1},\ldots, x_{n})\}\Big],\;\; n=0,1,2,\ldots
    \end{array}
\right.\label{eq:ADM1}
\end{equation}
Also the nonlinear term is decomposed as
 \begin{equation}
 Nx=\Big[1-\Big(\frac{dx}{dt}\Big)^2\Big]^{\frac{3}{2}}x=\sum_{n=0}^{\infty}A_{n}(x_{0},x_{1},\ldots, x_{n})
 \label{eq:N-1}
 \end{equation} 
where $\{A_{n}\}_{n=0}^{\infty}$ is the so-called Adomian polynomials sequence, the terms will be calculated according to \cite{Duan} and \cite{Duan1}. The first few polynomials are given by\\
\noindent $A_{0}(x_0)=x_{0}(1-x_{0}'^{2})^{\frac{3}{2}}, $\\
$A_{1}(x_0,x_1)=x_{1}(1-x_{0}'^{2})^{\frac{3}{2}},$\\
$A_{2}(x_0,x_1, x_2)=x_{2}(1-x_{0}'^{2})^{\frac{3}{2}}, $\\
$A_{3}(x_0,x_1, x_2, x_3)=x_{3}(1-x_{0}'^{2})^{\frac{3}{2}}, $\\
$A_{4}(x_0,x_1, x_2, x_3, x_4)=x_{4}(1-x_{0}'^{2})^{\frac{3}{2}}, $\\
$ \vdots $\\
$A_{m}(x_0,x_1,\ldots, x_m)=x_{m}(1-x_{0}'^{2})^{\frac{3}{2}} $ for every $m\geq 0$.\\
Now, recursively using (\ref{eq:ADM1}) with the Adomian polynomials given by the later sequence  $\{A_{n}\}_{n=0}^{\infty}$, we obtain, for a given initial velocity
$\beta$:\\
\begin{equation} \label{s-0}
x_{0}(t)=\beta t,
\end{equation}
\begin{equation} \label{s-1}
\begin{split}
x_{1}(t)&=-\mathcal{L}^{-1}\Big[\frac{1}{s^2}\mathcal{L}\{\beta(1-\beta^2)^{\frac{3}{2}}t\}\Big]=-\mathcal{L}^{-1}\Big[\frac{1}{s^4}\beta(1-\beta^2)^{\frac{3}{2}}\Big]\\ &=-\beta(1-\beta^2)^{\frac{3}{2}}\frac{t^3}{3!},
\end{split}
\end{equation}
\begin{equation} \label{s-2}
\begin{split}
x_{2}(t)&=-\mathcal{L}^{-1}\Big[\frac{1}{s^2}\mathcal{L}\{-\beta(1-\beta^2)^{3}\frac{t^3}{3!}\}\Big]=\mathcal{L}^{-1}\Big[\frac{1}{s^6}\beta(1-\beta^2)^{3}\Big]\\ &=\beta(1-\beta^2)^{3}\frac{t^5}{5!},
\end{split}
\end{equation}
\begin{equation} \label{s-3}
\begin{split}
x_{3}(t)&=-\mathcal{L}^{-1}\Big[\frac{1}{s^2}\mathcal{L}\{\beta(1-\beta^2)^{\frac{9}{2}}\frac{t^5}{5!}\}\Big]=-\mathcal{L}^{-1}\Big[\frac{1}{s^8}\beta(1-\beta^2)^{\frac{9}{2}}\Big]\\ &=-\beta(1-\beta^2)^{\frac{9}{2}}\frac{t^7}{7!},
\end{split}
\end{equation}
\begin{equation} \label{s-4}
\begin{split}
x_{4}(t)&=-\mathcal{L}^{-1}\Big[\frac{1}{s^2}\mathcal{L}\{-\beta(1-\beta^2)^{6}\frac{t^7}{7!}\}\Big]=\mathcal{L}^{-1}\Big[\frac{1}{s^{10}}\beta(1-\beta^2)^{6}\Big]\\ &=\beta(1-\beta^2)^{6}\frac{t^9}{9!},
\end{split}
\end{equation}
\begin{equation} \label{s-5}
\begin{split}
x_{5}(t)&=-\mathcal{L}^{-1}\Big[\frac{1}{s^2}\mathcal{L}\{\beta(1-\beta^2)^{\frac{15}{2}}\frac{t^9}{9!}\}\Big]=-\mathcal{L}^{-1}\Big[\frac{1}{s^{12}}\beta(1-\beta^2)^{\frac{15}{2}}\Big]\\ &=-\beta(1-\beta^2)^{\frac{15}{2}}\frac{t^{11}}{11!},
\end{split}
\end{equation}
$$\vdots . $$
In view of equations (\ref{s-0})-(\ref{s-5}), the series solution is
\begin{equation} \label{solser}
\begin{split}
x(t)&=\beta t-\beta(1-\beta^2)^{\frac{3}{2}}\frac{t^3}{3!}+\beta(1-\beta^2)^{3}\frac{t^5}{5!}-\beta(1-\beta^2)^{\frac{9}{2}}\frac{t^7}{7!}\\ & +\beta(1-\beta^2)^{6}\frac{t^9}{9!}-\beta(1-\beta^2)^{\frac{15}{2}}\frac{t^{11}}{11!}+\beta(1-\beta^2)^{9}\frac{t^{13}}{13!}\cdots\\
\end{split}
\end{equation}
$$=\beta\Big(t-(1-\beta^2)^{\frac{3}{2}}\frac{t^3}{3!}+(1-\beta^2)^{3}\frac{t^5}{5!}-(1-\beta^2)^{\frac{9}{2}}\frac{t^7}{7!}+(1-\beta^2)^{6}\frac{t^9}{9!}-+\cdots\Big)
$$
\begin{eqnarray}\label{SOL}
=\beta\sum_{n=0}^{\infty}\Big((1-\beta^2)^{\frac{3}{2}}\Big)^{n}Ç(-1)^{n}\frac{t^{2n+1}}{(2n+1)!}.
\end{eqnarray}
From (\ref{SOL}) we conclude that the solution of the equation  (\ref{Osc-1}), that is, the position of the relativistic harmonic oscillator is given by 
the series of power of functions with $0<\beta<1$
\begin{equation}\label{Ser}
x(t)=\beta\sum_{n=0}^{\infty}\Big((1-\beta^2)^{\frac{3}{2}}\Big)^{n}Ç(-1)^{n}\frac{t^{2n+1}}{(2n+1)!}.
\end{equation}
According to  \cite{Bart}, we easily see that the power series  (\ref{Ser}) converges in all $\mathbb{R}$ 
and it also converges uniformly in any compact subinterval of   $\mathbb{R}$.\\ 
Using the expressions obtained above for the solution of equation (\ref{Osc-1}), we will illustrate, with two examples, the
efectiveness of LADM to solve the nonlinear relativistic harmonic oscillator.

\section{\bf Application to the Relativistic Harmonic Oscillator}
\begin{figure}[h!]
	\begin{center}
		\includegraphics[width=110mm, height=65mm, scale=1.0]{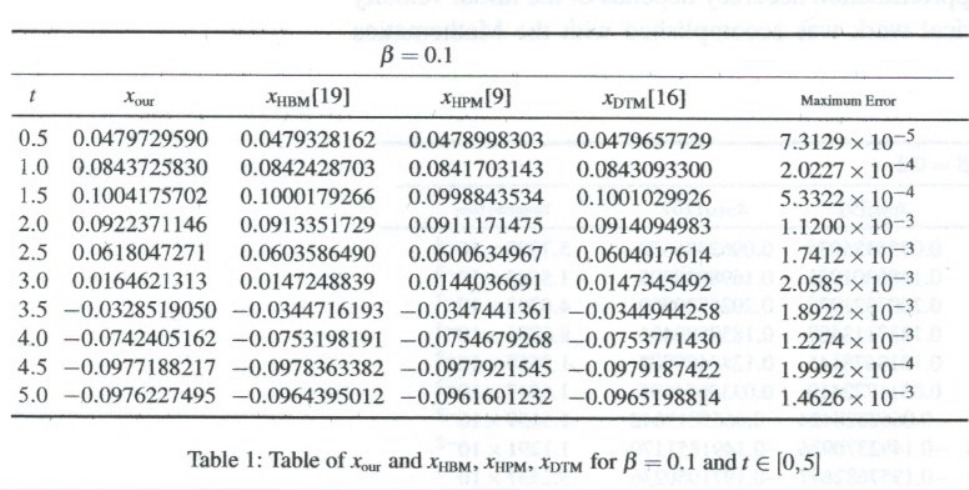}
	\end{center}
	\caption{Table for $\beta=0.1$ \label{tab1}}
\end{figure}

\begin{figure}[h!]
	\begin{center}
		\includegraphics[width=80mm, height=50mm, scale=1.0]{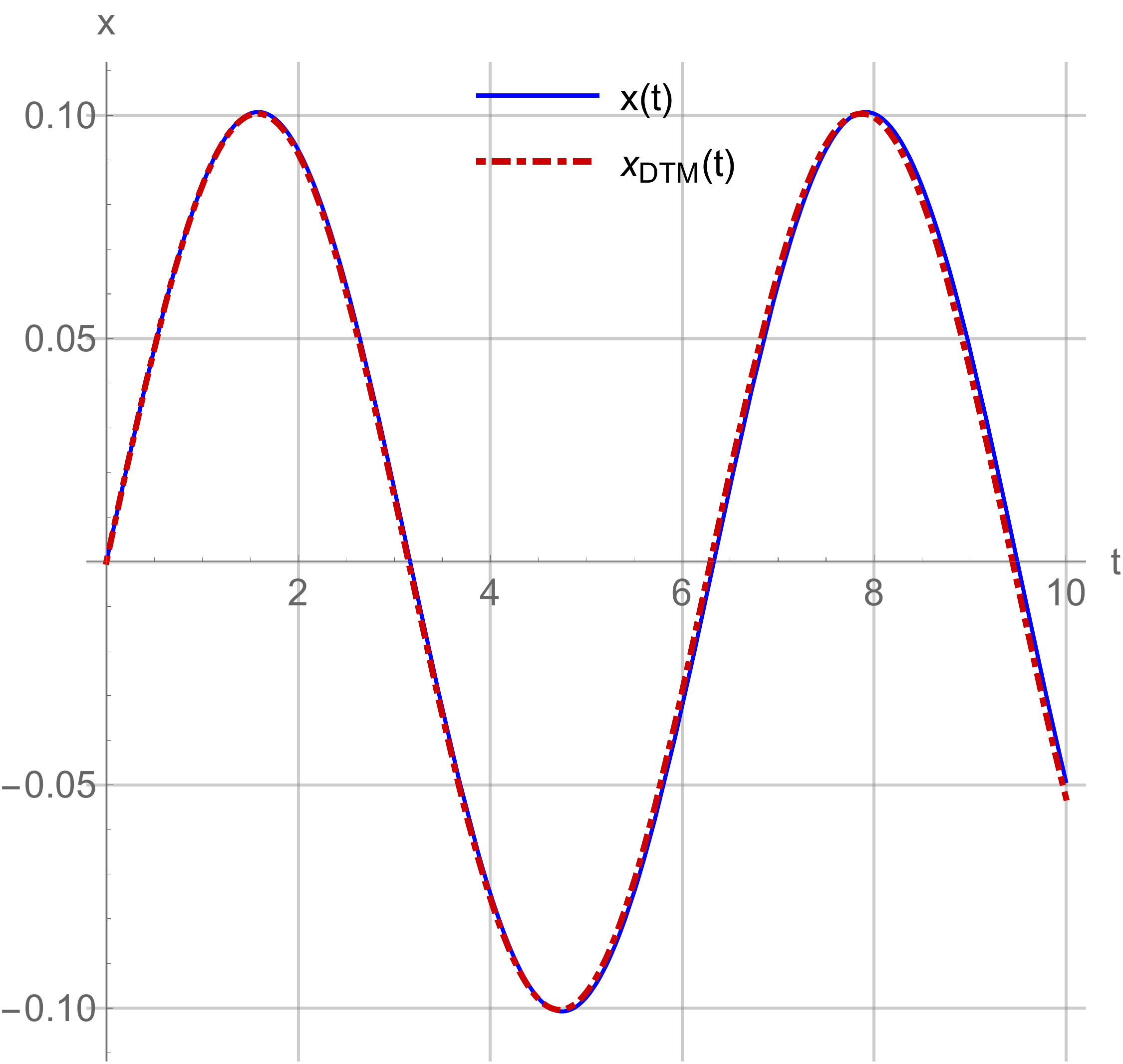}
	\end{center}
	\caption{Graph of the values of $x_{\mbox{\tiny our}}$ and $x_{\mbox{\tiny DTM}}$ for $\beta=0.1$ \label{fi1}}
\end{figure}

\begin{figure}[h!]
	\begin{center}
		\includegraphics[width=80mm, height=50mm, scale=1.0]{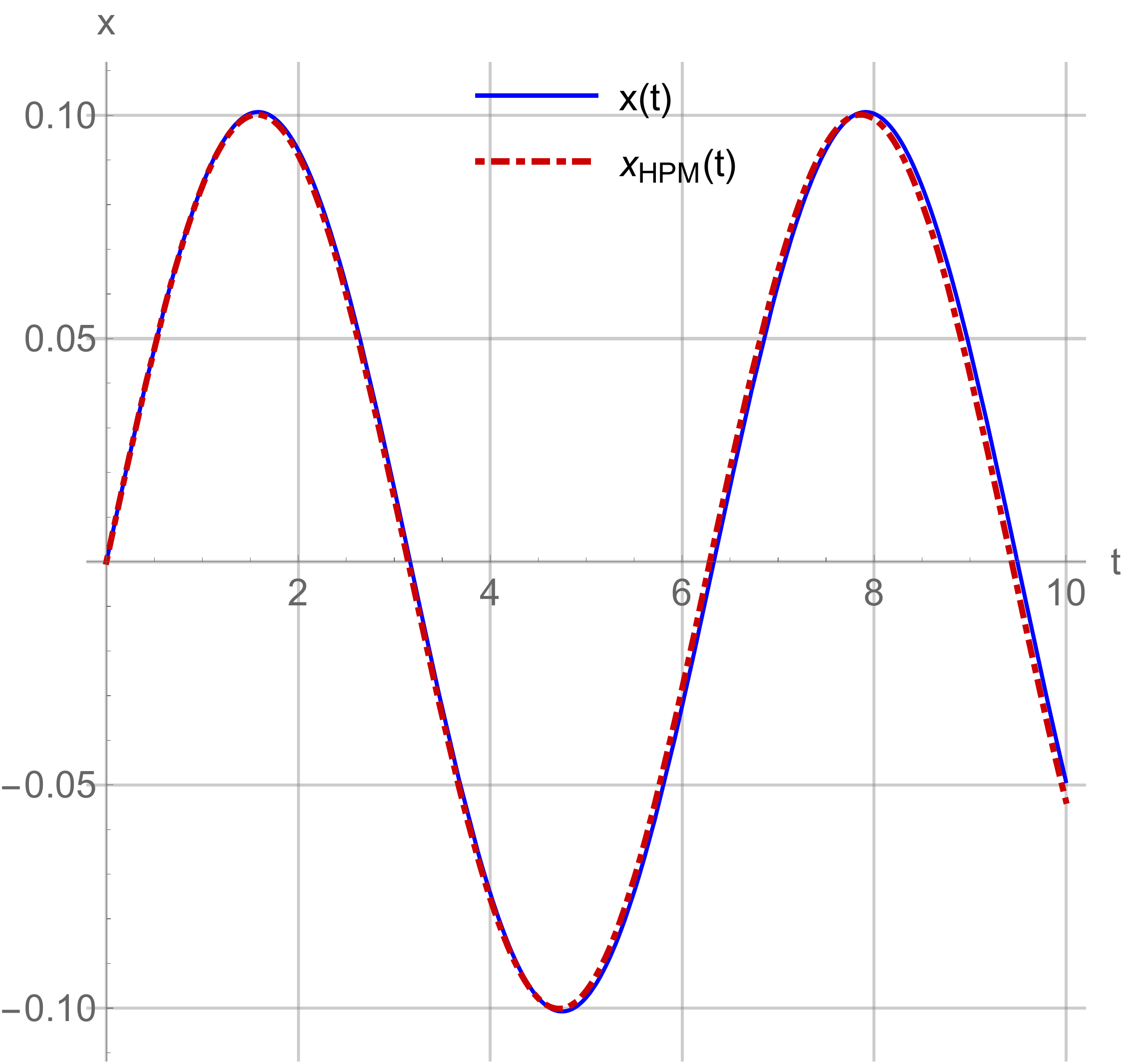}
	\end{center}
	\caption{Graph of the values of $x_{\mbox{\tiny our}}$ and $x_{\mbox{\tiny HPM}}$ for $\beta=0.1$ \label{fi2}}
\end{figure}

\begin{figure}[h!]
	\begin{center}
		\includegraphics[width=80mm, height=50mm, scale=1.0]{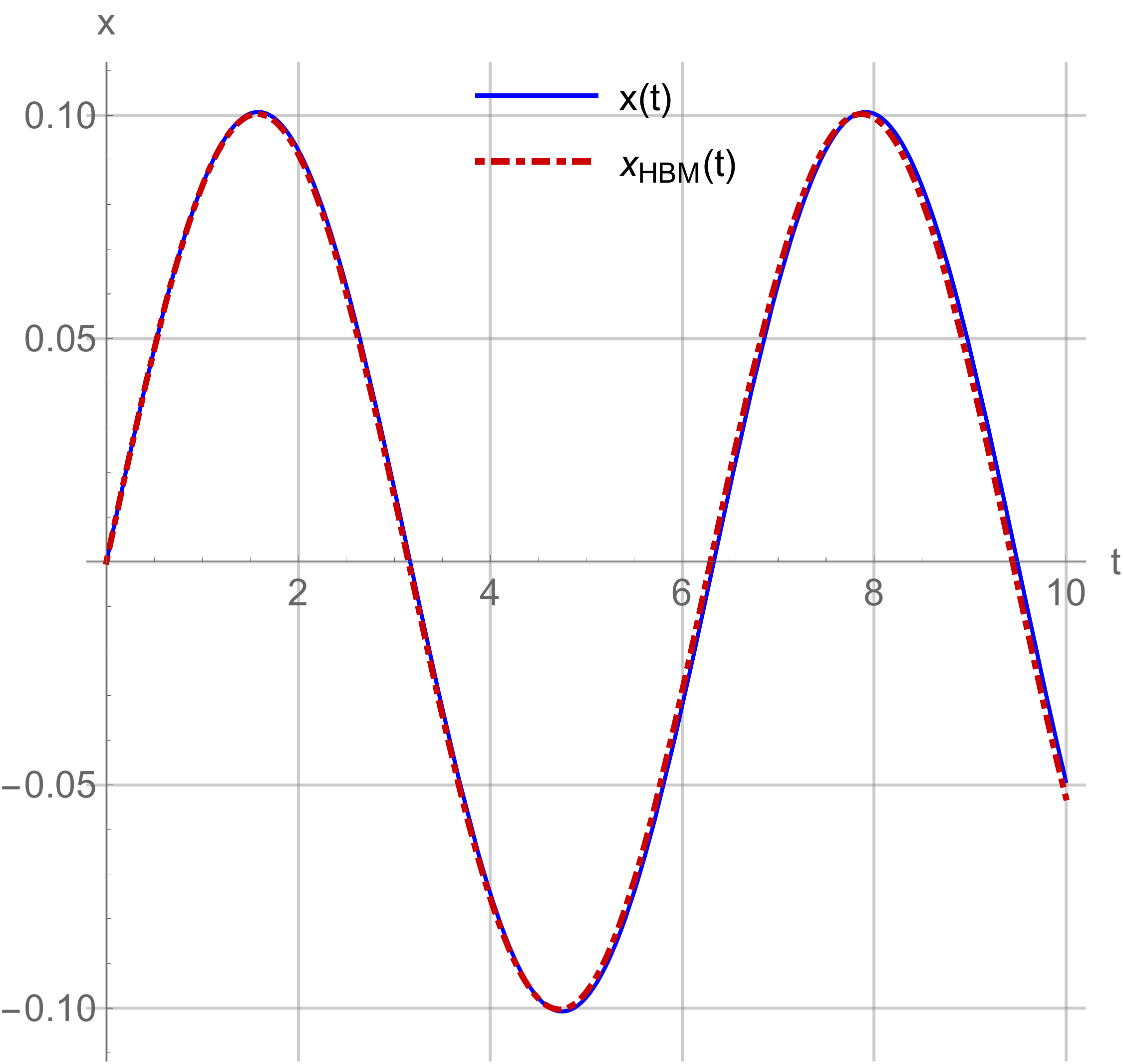}
	\end{center}
	\caption{Graph of the values of $x_{\mbox{\tiny our}}$ and $x_{\mbox{\tiny HBM}}$ for $\beta=0.1$ \label{fi3}}
\end{figure}

\noindent {\bf Example 1}\\
In this first example, we consider the particular case of (\ref{Osc-1}) such that $\beta =0.1$; 
this case was studied in \cite{Eb} via differential transformation method (DTM) and also in \cite{Bia} through the homotopy perturbation method (HPM). 
Good approximations were obtained in both works in comparison with the first known approximation solution of 
(\ref{Osc-1}) obtained in \cite{Mic1} by the  harmonic balance method (HBM).  We will use the formula 
 (\ref{Ser}) taking only the first fourteen terms (since the next one will be very small)
\begin{equation} \label{s1}
\begin{split}
x(t)&=0.1\sum_{n=0}^{13}(0.9850375)^{n}Ç(-1)^{n}\frac{t^{2n+1}}{(2n+1)!}=0.1t-0.0985037\frac{t^3}{3!}+0.0970299\frac{t^5}{5!}\\
& -0.095578\frac{t^7}{7!}+0.094148\frac{t^9}{9!}-0.0927393\frac{t^{11}}{11!}+\cdots
-0.0822027\frac{t^{27}}{27!}
\end{split}
\end{equation}
The approximations obtained for  $\beta=0.1$ v\'ia  DTM in \cite{Eb} by using HPM in \cite{Bia} are respectively:
\begin{equation}\label{n1}
x_{\mbox{\tiny DTM}}(t)=0.10033\sin(0.998t)-0.000047097\sin(2.997t)+0.00000008254\sin(4.841t)
\end{equation}

\begin{equation}\label{n2}
x_{\mbox{\tiny HPM}}(t)=0.10010\sin(0.999t)-0.00004689\sin(2.997t)+0.00000005062\sin(4.995t)
\end{equation}
Moreover, using $\beta=0.1$ in (\ref{hbm}) we find
\begin{equation}\label{n3}
x_{\mbox{\tiny HBM}}(t)=0.10025\sin(0.998t)-0.00004173\sin(2.996t)+0.00000004369\sin(4.944t)
\end{equation}
\noindent  The results obtained are shown in Table  \ref{tab1} in which the comparison with the ones  obtained in  \cite{Eb}, \cite{Bia} and \cite{Mic1} using  DTM, HPM and HBM respectively has been done.  We also display in figures  \ref{fi1}, \ref{fi2} and \ref{fi3}  this comparison.  
All the numerical work was accomplished with the Mathematica software package.

%\begin{figure}
%	\begin{center}
%		\includegraphics[width=85mm, height=50mm, scale=1.1]{G-1-1.pdf}\hspace{0.3in}
%		\includegraphics[width=85mm, height=50mm, scale=1.1]{G-1-2.pdf}\hspace{0.3in}
%		\includegraphics[width=85mm, height=50mm, scale=1.1]{G-1-3.pdf}%\hspace{0.3in}
%	\end{center}
%	\caption{Graph of the values of $x_{\mbox{\tiny our}}$ {\it versus} $x_{\mbox{\tiny DTM}}$, $x_{\mbox{\tiny HPM}}$ and $x_{\mbox{\tiny HBM}}$ for $\beta=0.1$ \label{fi1}}
%\end{figure}

 \noindent {\bf Example 2}\\
 In this second example, we consider the particular case of (\ref{Osc-1}) such that $\beta =0.2$; this case was studied in \cite{Eb} via DTM and also in  \cite{Bia}  using 
 HPM.  Once again, in both works, good approximations were found in comparison with the first obtained in  (\ref{Osc-1}) and the one obtained in \cite{Mic1} by HBM. 
 As before, using the formula  (\ref{Ser}) taking the first fourteen terms we obtain 
 \begin{equation} \label{s1}
 \begin{split}
 x(t)&=0.2\sum_{n=0}^{13}(0.940604)^{n}Ç(-1)^{n}\frac{t^{2n+1}}{(2n+1)!}=0.2t-0.1881208\frac{t^3}{3!}+0.1769472\frac{t^5}{5!}\\
 & -0.1664372\frac{t^7}{7!}+0.1565515\frac{t^9}{9!}-0.1472253\frac{t^{11}}{11!}+\cdots
 -0.0902233\frac{t^{27}}{27!}
 \end{split}
 \end{equation}
 The approximations obtained in the case of $\beta=0.2$ v\'ia  DTM in \cite{Eb} through HPM in \cite{Bia} are respectively:
 
 \begin{equation}\label{p1}
 x_{\mbox{\tiny DTM}}(t)=0.203\sin(0.992t)-0.0003695\sin(3.051t)+0.000009257\sin(4.29t) 
  \end{equation}

  \begin{equation}\label{p2}
 x_{\mbox{\tiny HPM}}(t)=0.201\sin(0.995t)-0.0003768\sin(2.985t)+0.000001652\sin(4.974t) 
  \end{equation}
And also using    $\beta=0.2$ in (\ref{hbm}) we obtain
 
  \begin{equation}\label{p3}
x_{\mbox{\tiny HBM}}(t)=0.202\sin(0.995t)-0.0003354\sin(2.985t)+0.000001508\sin(4.974t)  
  \end{equation}

\begin{figure}[h!]
	\begin{center}
		\includegraphics[width=110mm, height=65mm, scale=1.0]{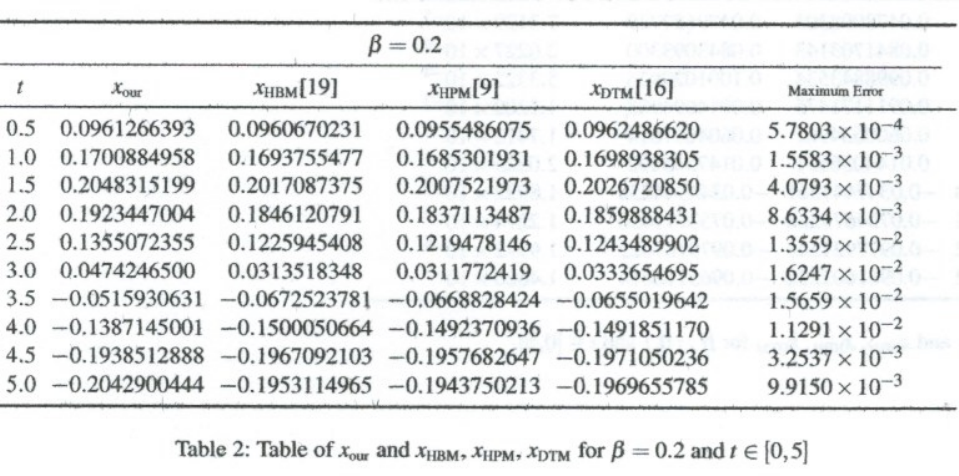}
	\end{center}
	\caption{Table for $\beta=0.2$ \label{tab2}}
\end{figure}

\begin{figure}[h!]
	\begin{center}
		\includegraphics[width=80mm, height=50mm, scale=1.0]{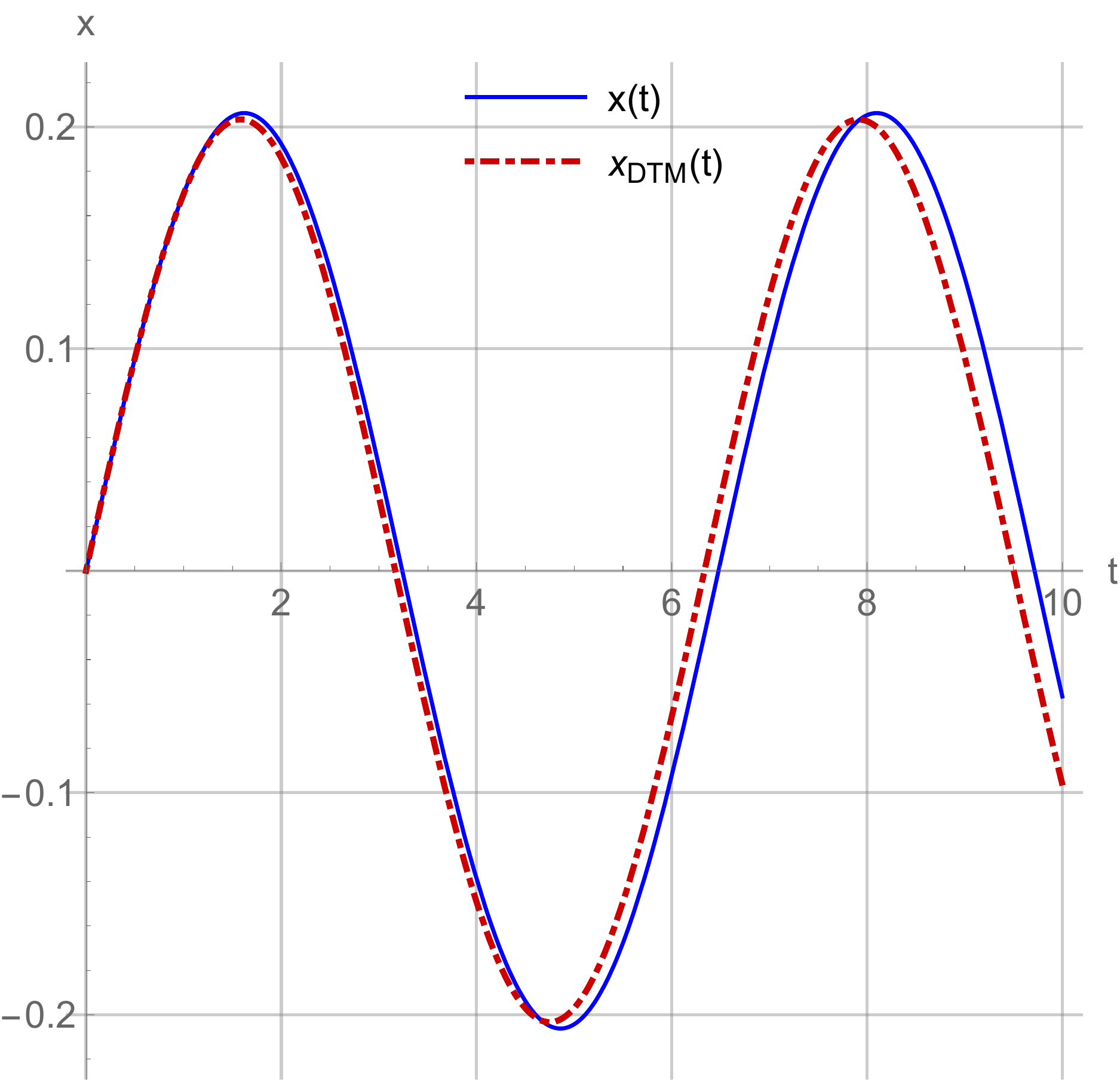}
	\end{center}
	\caption{Graph of the values of $x_{\mbox{\tiny our}}$ and $x_{\mbox{\tiny DTM}}$ for $\beta=0.2$ \label{fi4}}
\end{figure}

\begin{figure}[h!]
	\begin{center}
		\includegraphics[width=80mm, height=50mm, scale=1.0]{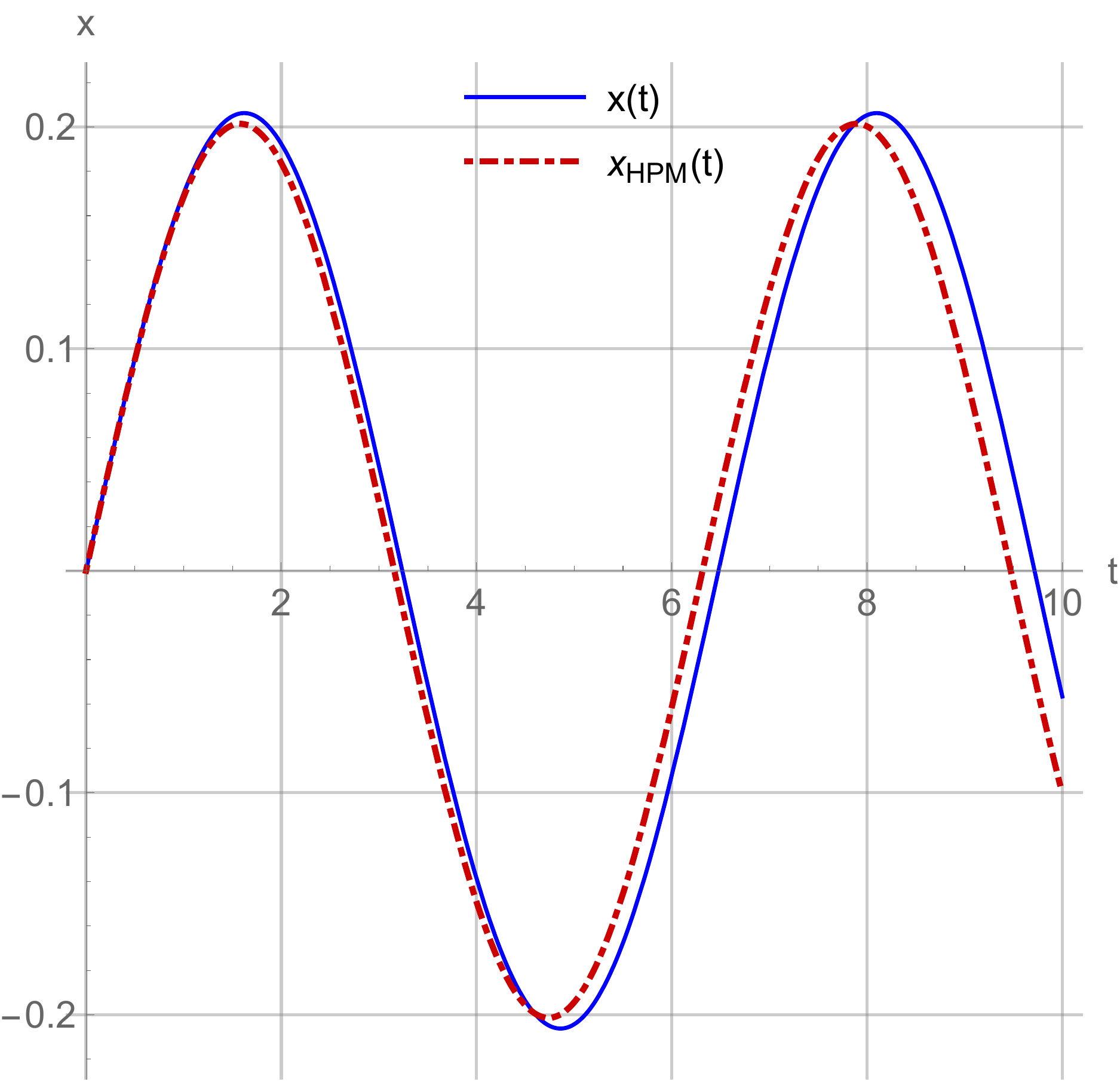}
	\end{center}
	\caption{Graph of the values of $x_{\mbox{\tiny our}}$ and $x_{\mbox{\tiny HPM}}$ for $\beta=0.2$ \label{fi5}}
\end{figure}

\begin{figure}[h!]
	\begin{center}
		\includegraphics[width=80mm, height=50mm, scale=1.0]{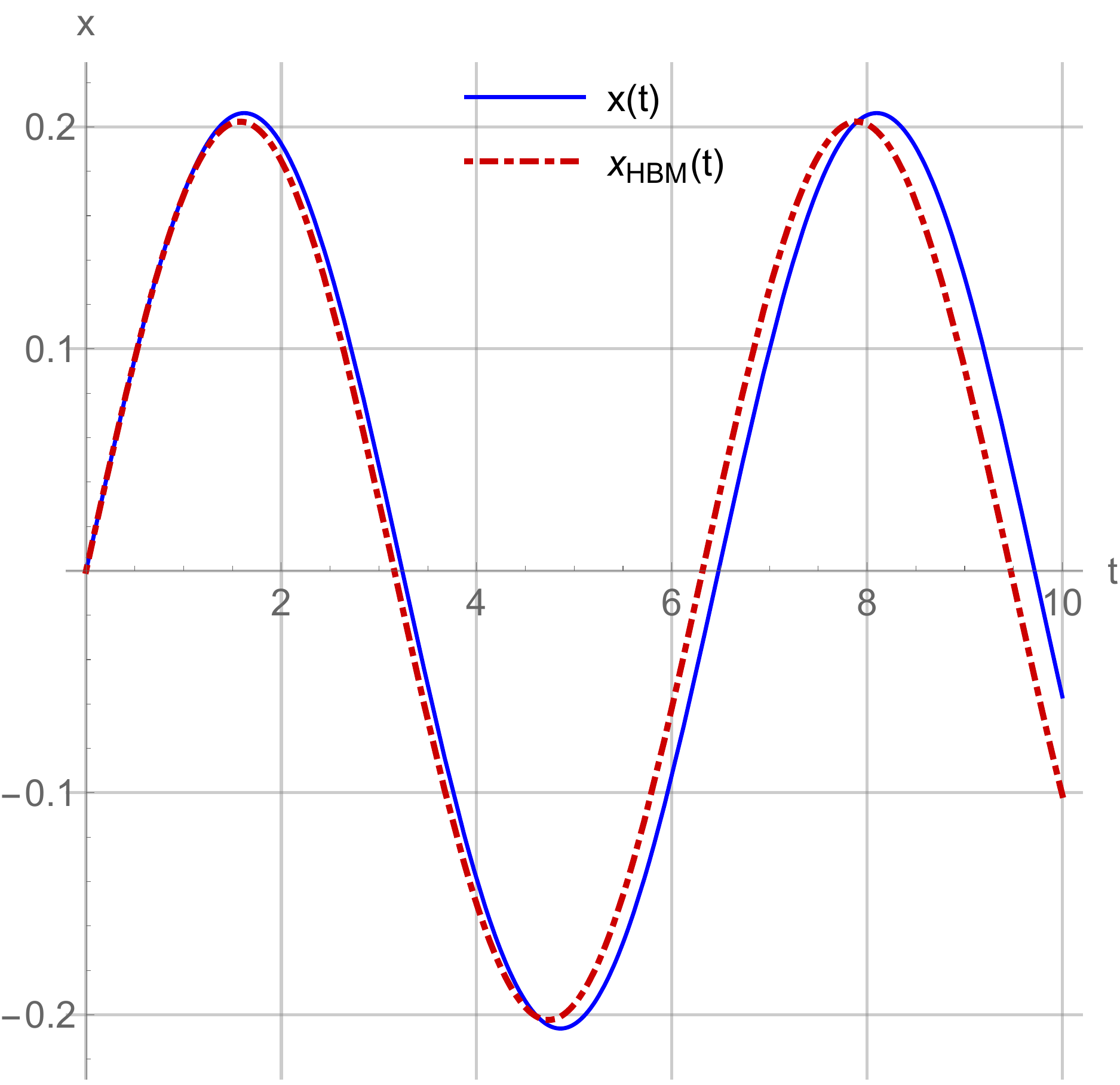}
\end{center}
	\caption{Graph of the values of $x_{\mbox{\tiny our}}$ and $x_{\mbox{\tiny HBM}}$ for $\beta=0.2$ \label{fi6}}
\end{figure}

 \noindent Comparison of our results with the ones obtained in \cite{Eb}, \cite{Bia} and \cite{Mic1} using  DTM, HPM and HBM are showed in Table \ref{tab2} and displayed 
 in figures \ref{fi4}, \ref{fi5} and \ref{fi6}. In this example we can also see that the approximation accuracy  depends of the initial  velocity of the oscillator.
 All the numerical work was accomplished with the Mathematica software package.

%\begin{figure}
%	\begin{center}
%		\includegraphics[width=85mm, height=50mm, scale=1.1]{G-2-1.pdf}\hspace{0.3in}
%		\includegraphics[width=85mm, height=50mm, scale=1.1]{G-2-2.pdf}\hspace{0.3in}
%		\includegraphics[width=85mm, height=50mm, scale=1.1]{G-2-3.pdf}%\hspace{0.3in}
%	\end{center}
%	\caption{Graph of the values of $x_{\mbox{\tiny our}}$ {\it versus} $x_{\mbox{\tiny DTM}}$, %$x_{\mbox{\tiny HPM}}$ and $x_{\mbox{\tiny HBM}}$ for $\beta=0.2$ \label{fi2}}
%\end{figure}

\noindent As we seen from the last examples, the solutions  we have obtained are periodic functions and the amplitude depends of the initial velocity 
as found by the author in  \cite{Mic1}. The main difference of our results with the reported ones is that the final series is uniformly convergent in any compact subset 
of the real line and therefore we can obtain the results with the required accuracy.

\section{Summary and Conclusions}

\noindent Far as we know there is no exact solutions to the equation of motion for
a relativistic harmonic oscillator. In this work, we have obtained the solution of the problem without 
the {\it a priori} assumption that the solutions are periodic functions; 
the solution that we have obtained is a series of powers of functions which uniformly converge on compact subsets of $\mathbb{R}$, never before reported. 
The problem of find the limit function of the series solution is an open question that we are currently work.\\
\noindent In order to show the accuracy and efficiency of our method,  we have solved two examples and comparing our results with the ones obtained with three different methods \cite{Mic1,Eb,Bia}. Our results show that LADM produces highly accurate solutions in complicated nonlinear problems. We therefore,  conclude that the Laplace-Adomian decomposition method is a notable  non-sophisticated powerful tool that produces high quality approximate solutions for nonlinear ordinary differential equations using simple calculations and that reaches convergence with only  a few terms.  Finally, the Laplace-Adomian decomposition method would be a powerful mathematical tool for solving other nonlinear differential equations 
related with mathematical physics models.  All the numerical work and the graphics was accomplished with the Mathematica software package.

%\section*{Acknowledgments}
%We would like to thank anonymous referees for their constructive comments and suggestions that %helped to improve the paper.

\end{document}